\documentclass{amsart}
\usepackage{amsmath,amssymb,amsthm,mathrsfs}

\newcommand{\BC}{{\mathbb C}}\newcommand{\BD}{{\mathbb D}}

\newcommand{\sD}{{\mathcal D}}
\newcommand{\sF}{{\mathcal F}}
\newcommand{\sG}{{\mathcal G}}\newcommand{\sH}{{\mathcal H}}

\newcommand{\sK}{{\mathcal K}}
\newcommand{\sM}{{\mathcal M}}

\newcommand{\sU}{{\mathcal U}}\newcommand{\sV}{{\mathcal V}}
\newcommand{\sW}{{\mathcal W}}\newcommand{\sX}{{\mathcal X}}
\newcommand{\sY}{{\mathcal Y}}

\newcommand{\eH}{{\mathbf H}}

\newcommand{\eS}{{\mathbf S}}

\newcommand{\cE}{{\mathscr E}}

\newcommand{\Ga}{\Gamma}
\newcommand{\de}{\delta}\newcommand{\De}{\Delta}

\newcommand{\la}{\lambda}

\newcommand{\Up}{\Upsilon}

\newcommand{\om}{\omega}

\newcommand{\kr}{\textup{Ker\,}}

\newcommand{\mat}[2]{\ensuremath{\left[\begin{array}{#1}
#2
\end{array} \right]}}
\newcommand{\ov}[1]{{\overline{#1}}}

\newcommand{\tu}[1]{\textup{#1}}

\newcommand{\ands}{\quad\mbox{and}\quad}
\newcommand{\ors}{\quad\mbox{or}\quad}

\newcommand{\LDS}{\{A,T',U',R,Q\}}

\newtheorem{theorem}{Theorem}[section]
\newtheorem{corollary}[theorem]{Corollary}
\newtheorem{lemma}[theorem]{Lemma}
\newtheorem{proposition}[theorem]{Proposition}

\newtheorem{example}[theorem]{Example}

\begin{document}

\title{Relaxed commutant lifting: existence of a unique solution}

\author{S. ter Horst}
\address{Department of Mathematics, Virginia Tech, Blacksburg, VA 24061-0123, USA}
\email{terhorst@math.vt.edu}

\date{}

\subjclass{47A20, 47A56, 47A57}
\keywords{commutant lifting, unique solution, linear fractional maps}

\begin{abstract}
In this paper we present necessary and sufficient conditions for the existence of a unique
solution to the relaxed commutant lifting problem.
The obtained conditions are more complicated than those for the classical commutant lifting
setting, and earlier obtained sufficient conditions turn out not to be necessary conditions.
It is also shown that these conditions simplify in certain special cases.

\end{abstract}

\maketitle

\section*{Introduction}\label{sec:intro}

In this paper we resolve an open problem concerning relaxed commutant lifting.
The relaxed commutant lifting problem was introduced in \cite{FFK02}, extending the
classical commutant lifting theory of B. Sz.-Nagy and C. Foias \cite{NF68} and its
predecessor by D. Sarason \cite{S67}. As mentioned in \cite{FFK02}, relaxed commutant lifting
covers the generalizations of the commutant lifting setting by Treil and Volberg \cite{TV94}
and the weighted version in \cite{BFF99}. In \cite{FFK02} it was shown that a particular (central)
solution always exists. Representations of all solutions have been given in
\cite{FtHK06a,LT06,FtHK06b,tH}.
In the present paper we provide necessary and sufficient conditions under which the
so-called central solution is the only solution to the relaxed commutant lifting problem.

The data for the relaxed commutant lifting problem is a set $\LDS$ consisting of five
Hilbert space operators: the operator $A$ is a contraction mapping $\sH$ into ${\sH}'$, the
operator $U'$ on $\sK'$ is a minimal isometric lifting of the contraction $T'$ on $\sH'$, i.e.,
$U'$ is an isometry on $\sK'$, with $\sH'\subset\sK'$ being cyclic for $U'$ and
$\Pi_{\sH'}U'=T'\Pi_{\sH'}$, and $R$ and $Q$ are operators from $\sH_0$ to $\sH$, satisfying
\begin{equation}\label{intertw}
T'AR=AQ\ands R^*R\leq Q^*Q.
\end{equation}
Here we follow the convention that for a subspace $\sV$ of a Hilbert space $\sW$ the symbol
$\Pi_\sV$ stands for the orthogonal projection from $\sW$ onto $\sV$ viewed as an operator
from $\sW$ onto $\sV$.
Given this data set the {\em relaxed commutant lifting problem} is to find a (all) contraction(s)
$B$ from $\sH$ to $\sK'$ such that
\begin{equation}\label{rclt}
\Pi_{\sH'}B=A\ands U'BR=BQ.
\end{equation}
A contraction $B$ from $\sH$ into $\sK'$ that satisfies (\ref{rclt}) is called a
{\em contractive interpolant for $\LDS$}.

An essential role in the solution to the relaxed commutant lifting problem is played by the
so-called {\em underlying contraction} $\om$ (see \cite{FFK02}), which is defined by
\begin{equation}\label{omdef}
\om:\sF=\ov{D_AQ\sH_0}\to\mat{c}{\sD_{T'}\\\sD_A},\quad \om
D_AQ=\mat{c}{D_{T'}AR\\D_AR\\}.
\end{equation}
Here, as usual, given a contraction $N$, we write $D_N$ for the defect operator and $\sD_N$
for the defect space of $N$, that is, $D_N$ is the positive square root of $I-N^*N$ and
$\sD_N$ is the closure of the range of $D_N$.
The first component of $\om$, mapping $\sF$ into $\sD_{T'}$, is denoted by $\om_1$, and
the second component of $\om$, mapping $\sF$ into $\sD_{A}$, by $\om_2$.


The classical commutant lifting setting appears when $R$ is the identity operator on $\sH$, and
thus $\sH_0=\sH$, and $Q$ is an isometry.
It is well known that in the classical commutant lifting setting there exists a unique contractive
interpolant if and only if $\ov{\sD_A Q\sH_0}=\sD_A$ or $\om$ is unitary, cf., the first paragraph
after Theorem VI.2.1 in \cite{FFGK98}.

For relaxed commutant lifting these conditions are sufficient, but not necessary;
see Example \ref{ex:counter} below. The following theorem provides conditions
for the relaxed commutant lifting setting that are both necessary and sufficient.

\begin{theorem}\label{mainth1}
There exists a unique contractive interpolant for the relaxed commutant lifting problem
with data set $\LDS$ if and only if the underlying contraction $\om$ satisfies one of
the following conditions:
\begin{itemize}

\item[\tu{(i)}] $\ov{\sD_A Q\sH_0}=\sD_A$;

\item[\tu{(ii)}] the operator $\om_1(\Pi_\sF\om_2)^n$ is a
co-isometry for $n=0,1,2,\ldots$.

\end{itemize}
\end{theorem}

Condition (i) was already obtained as a sufficient condition; see Theorem 3.1 in \cite{FFK02}.
Two other sufficient conditions were also obtained earlier; both 
are covered by (ii). Namely,
the trivial condition: $T'$ is an isometry ($\sD_{T'}=\{0\}$) (Proposition 6.6 in \cite{FtHK06a}), and the condition that
$\om$ is a co-isometry (Section 4.1 in \cite{tH07}).

To see that this theorem covers the classical result it suffices to show that
condition (ii) in the classical commutant lifting setting implies that $\om$ is unitary;
the sufficiency of the conditions follows from the last remark in the previous paragraph.
Note that in the classical case $R^*R=Q^*Q$, which implies that the underlying contraction $\om$
is an isometry; see \cite{FFK02}. Moreover, since $R$ is the identity operator, the range of $\om_2$
is dense in $\sD_A$. So if condition (ii) holds, then in particular $\om_1$ is a co-isometry,
and, as $\om$ is an isometry and $\ov{\om_2\sF}=\sD_A$, the contraction $\om_2$ is also a co-isometry.
This implies that $\om$ is a co-isometry, and thus unitary.

In order to prove Theorem \ref{mainth1} it is more convenient to work with the equivalent
interpolation problem of \cite{FtHK3}.
To state the latter problem we first need some preliminaries.
Let $\sU$ and $\sY$ be Hilbert spaces. We
write $\eH^2(\sU,\sY)$ for the set of all operator-valued functions $H$ on the unit disc
$\BD$ of $\BC$ whose values are operators from $\sU$ into $\sY$, and such that the formula
\begin{equation}\label{GaH}
(\Ga_H u)(\la)=H(\la)u\qquad(u\in\sU,\la\in\BD)
\end{equation}
defines an operator $\Ga_H$ from $\sU$ into the Hardy space
$H^2(\sY)$. The set $\eH^2(\sU,\sY)$ is a Banach space under the
norm $\|H\|=\|\Ga_H\|$, where $\|\Ga_H\|$ denotes the operator norm
of $\Ga_H$. We write $\eH^2_\tu{ball}(\sU,\sY)$ for the unit ball of
$\eH^2(\sU,\sY)$.

With these definitions the interpolation problem of \cite{FtHK3} can be formulated as follows:
given a contraction
\begin{equation}\label{givenom}
\om=\mat{c}{\om_1\\\om_2}:\sF\to\mat{c}{\sY\\\sU},
\quad\mbox{where }\sF\text{ is a subspace of }\sU,
\end{equation}
describe a (all) function(s) $H$ in $\eH^2_\tu{ball}(\sU,\sY)$ such that
\begin{equation}\label{intercon}
\om_1 +\la H(\la) \om_2 =H(\la)|_\sF\quad(\la\in\BD).
\end{equation}
A function $H$ in $\eH^2_\tu{ball}(\sU,\sY)$ satisfying (\ref{intercon}) is called a
{\em solution to the $H^2$ interpolation problem defined by (the contraction) $\om$},
or just a solution when the contraction $\om$ in question is clear from the context.

Notice that the underlying contraction $\om$ in (\ref{omdef}) is of the form (\ref{givenom}).
By making an appropriate choice for the data set, each contraction $\om$ of the
form (\ref{givenom}) appears as the underlying contraction for the corresponding  relaxed
commutant lifting problem.
Moreover, given a data set $\LDS$ with underlying contraction $\om$, there is a one-to-one
map between the solutions $H$ of the $H^2$ interpolation problem defined by $\om$
and the contractive interpolants $B$ for $\LDS$; this map is given by
\begin{equation*}
H\mapsto B=\mat{c}{A\\\Ga_H D_A},
\end{equation*}
provided the minimal isometric lifting $U'$ is the Sz.-Nagy-Sch\"affer isometric lifting
on $\sH'\oplus H^2(\sD_{T'})$,
which we may assume without loss of generality. In this sense the two problems are
equivalent; see \cite{FtHK3} for details.


The following theorem is the analog of Theorem \ref{mainth1} for the $H^2$ interpolation setting.

\begin{theorem}\label{mainth2}
Let $\om$ be a contraction of the form $(\ref{givenom})$. Then there exists a unique
solution to the $H^2$ interpolation problem defined by $\om$
if and only if one of the following conditions is met:
\begin{itemize}

\item[\tu{(i)}] $\sF=\sU$;

\item[\tu{(ii)}] the operator $\om_1(\Pi_\sF\om_2)^n$ is a
co-isometry for $n=0,1,2,\ldots$.

\end{itemize}
\end{theorem}

We prove Theorem \ref{mainth2} in Section \ref{sec:UniqueSolution}.
Theorem \ref{mainth1} follows directly from Theorem \ref{mainth2}
and the equivalence between the two problems, as explained above.

A solution to the $H^2$ interpolation problem defined by $\om$ always exists.
Indeed, it is not difficult to verify that the function
$H_c$ on $\BD$ given by
\begin{equation}\label{centralCI}
H_c(\la)=\om_1\Pi_\sF(I_{\sU}-\la\om_2\Pi_\sF)^{-1}\qquad(\la\in\BD)
\end{equation}
is in $\eH^2_\tu{ball}(\sU,\sY)$ and satisfies (\ref{intercon}).
This particular solution $H_c$ is referred to as the {\em central solution}, and corresponds to the
central solution for the relaxed commutant lifting problem obtained in \cite{FFK02}.

The operators $\om_1(\Pi_\sF\om_2)^n$ in condition (ii) are closely related to the Taylor
coefficients of the central solution $H_c$. Indeed, the $n^\tu{th}$ Taylor
coefficient of $H_c$ is equal to $\om_1\Pi_\sF(\om_2\Pi_\sF)^n=\om_1(\Pi_\sF\om_2)^n\Pi_\sF$,
and thus condition (ii) holds if and only if the Taylor coefficients of $H_c$ are co-isometric.
Now let $\Ga_{H_c}$ be the contraction from $\sU$ into $H^2(\sY)$ associated with $H_c$ via
(\ref{GaH}). It then follows that condition (ii) in Theorem \ref{mainth2} is equivalent to
\begin{itemize}

\item[$\tu{(ii)}'$] the operator $\Ga_{H_c}$ is a co-isometry.

\end{itemize}
In particular, we have the following corollary.

\begin{corollary}\label{cor:coisosolution}
Let $\om$ be a contraction of the form $(\ref{givenom})$, and assume that $\sF\not=\sU$.
Then $H_c$ is the only solution if and only if $\Ga_{H_c}$ is a co-isometry.
\end{corollary}

This paper consists of three sections, not counting the present introduction. The first
section contains preliminary material on a Redheffer representation of all solutions.
In Section 2 we prove Theorem \ref{mainth2}, give an example and show how condition (ii)
in Theorem \ref{mainth2} simplifies when additional conditions are imposed upon $\om$.
In the final section we revisit the original relaxed commutant lifting setting and specify
Theorem \ref{mainth1} for some metric constrained interpolation problems.

We conclude the introduction with some words on the notation and terminology used in this paper.
Throughout calligraphic letters denote Hilbert spaces. We write $\dim(\sU)$ for the dimension of $\sU$.
The Hilbert space direct sum of $\sU$ and
$\sY$ is denoted by
\begin{equation*}
\sU\oplus\sY\ors\mat{c}{\sU\\\sY}.
\end{equation*}
By definition, a {\em subspace} is a closed linear manifold, and the closure of a linear subset
$\cE$ of $\sU$ is denoted by $\ov{\cE}$. If $\sM$ is a subspace of $\sU$, then $\sU\ominus\sM$
stands for the orthogonal complement of $\sM$ in $\sU$. We write $\sM_1\vee\sM_2$ for the
closure of the linear span of the subspaces $\sM_1$ and $\sM_2$ of $\sU$.
The term {\em operator} stands for a bounded linear transformation acting between Hilbert spaces.
We say that a contraction $A$ is a strict contraction if
$\|A\|<1$, or equivalently, when the defect operator $D_A$ is invertible. An operator $C$ from $\sU$ into $\sY$
is said to be {\em left invertible} whenever there exists an operator $D$ from $\sY$ to $\sU$ with
$DC=I_\sU$ (the identity operator on $\sU$). In this case $D$ is called a {\em left inverse} of $C$.
We say that a subspace $\sM$ of $\sU$ is {\em cyclic} for an operator $T$ on $\sU$ if the smallest
subspace of $\sU$ that contains $T^n\sM$ for each nonnegative integer $n$ is $\sU$ itself.
Finally, a {\em Schur class function} is an operator-valued function on the unit disc $\BD$
whose values are contractions. The set of Schur class functions whose values are operators
from $\sU$ into $\sY$ is denoted by $\eS(\sU,\sY)$, and is referred to as the {\em Schur class
associated with $\sU$ and $\sY$}.

\section{A Redheffer representation of all solutions}\label{sec:Redheffer}

Let $\om$ be a contraction of the form (\ref{givenom}). In this section we review some facts
concerning the Redheffer representation of \cite[Section 5.1]{tH07} that describes all
solutions to the $H^2$ interpolation problem defined by $\om$.

Set $\sG=\sU\ominus\sF$, and define operator-valued functions $\Phi_{1,1}$, $\Phi_{1,2}$,
$\Phi_{2,1}$ and $\Phi_{2,2}$ on $\BD$ by
\begin{equation}\label{REDcoefs}
\begin{array}{rl}
&\Phi_{1,1}(\la)=\la\Pi_\sG(I_\sU-\la\om_2\Pi_\sF)^{-1}
\Pi_\sU D_{\om^*},\\[.1cm]
&\Phi_{1,2}(\la)=\Pi_\sG(I_\sU-\la\om_2\Pi_\sF)^{-1},\\[.1cm]
&\Phi_{2,1}(\la)=\Pi_\sY D_{\om^*}+\la\om_1\Pi_\sF
(I_\sU-\la\om_2\Pi_\sF)^{-1}\Pi_\sU D_{\om^*},\\[.1cm]
&\Phi_{2,2}(\la)=\om_1\Pi_\sF(I_\sU-\la\om_2\Pi_\sF)^{-1}.
\end{array}\quad(\la\in\BD)
\end{equation}
Using the equivalence between relaxed commutant lifting and the $H^2$ interpolation
problem defined by $\om$, as described in \cite{FtHK3}, Theorem 5.1.1 from \cite{tH07} provides
the following Redheffer representation of all solutions.

\begin{theorem}\label{th:Redheffer}
Let $\om$ be a contraction of the form \tu{(\ref{givenom})}, and let $\Phi_{1,1}$, $\Phi_{1,2}$,
$\Phi_{2,1}$ and $\Phi_{2,2}$ be the functions given by \tu{(\ref{REDcoefs})}. Then
\begin{equation*}
\Phi_{1,1}\in\eS(\sD_{\om^*},\sG),\quad\Phi_{1,2}\in\eH^2_\tu{ball}(\sU,\sG)
,\quad\Phi_{2,1}\in\eS(\sD_{\om^*},\sY),\quad\Phi_{2,2}\in\eH^2_\tu{ball}(\sU,\sY),
\end{equation*}
and for any $V\in\eS(\sG,\sD_{\om^*})$ the function $H$ given by
\begin{equation}\label{REDsols}
H(\la)=
\Phi_{2,2}(\la)+\Phi_{2,1}(\la)V(\la)
(I-\Phi_{1,1}(\la)V(\la))^{-1}\Phi_{1,2}(\la)
\quad(\la\in\BD)
\end{equation}
is a solution to the $H^2$ interpolation problem defined by $\om$.
Moreover, all solutions are obtained in this way.
\end{theorem}

A representation of the form (\ref{REDsols}) is referred to as a {\em Redheffer representation},
and the functions $\Phi_{1,1}$, $\Phi_{1,2}$, $\Phi_{2,1}$ and $\Phi_{2,2}$ as the corresponding
{\em Redheffer coefficients}.

Since $\Phi_{1,1}$ is a Schur class function with $\Phi_{1,1}(0)=0$, it follows that
$\Phi_{1,1}(\la)=\la\Up(\la)$ for some Schur class function $\Up$ in $\eS(\sD_{\om^*},\sG)$
and for each $\la\in\BD$; see Lemma 2.4.1 in \cite{tH07}. Thus $\Phi_{1,1}(\la)$ is a strict
contraction for each $\la\in\BD$ and the inverse in (\ref{REDsols}) is properly defined.

By taking $V\in\eS(\sG,\sD_{\om^*})$ to be the zero function, we see that the function
$\Phi_{2,2}$ is a solution. In fact, $\Phi_{2,2}$ is precisely the central solution $H_c$ in
(\ref{centralCI}).

The map $V\mapsto H$ given by the Redheffer representation (\ref{REDsols}) is in general not
one-to-one. It can happen that different $V$'s yield the same solution $H$. However, this
non-uniqueness in the representation can be made explicit; see \cite{FtHK06b,tH07} for details.

Since $\Phi_{1,1}$ and $\Phi_{2,1}$ are Schur class functions, they define, in the usual way,
contractive
multiplication operators $M_{\Phi_{1,1}}$ from $H^2(\sD_{\om^*})$ into $H^2(\sG)$ and
$M_{\Phi_{2,1}}$ from $H^2(\sD_{\om^*})$ into $H^2(\sY)$, respectively. Moreover,
Theorem \ref{th:Redheffer} implies that we can define contractions $\Ga_{\Phi_{1,2}}$ from $\sU$
into $H^2(\sG)$ and $\Ga_{\Phi_{2,2}}$ from $\sU$ into $H^2(\sY)$ associated with $\Phi_{1,2}$,
respectively $\Phi_{2,2}$, via (\ref{GaH}).

The following theorem provides an improved version of Proposition 5.1.6 in \cite{tH07}.

\begin{theorem}\label{th:REDcoefmat}
Let $\om$ be a contraction of the form \tu{(\ref{givenom})}, let $M_{\Phi_{1,1}}$ and
$M_{\Phi_{2,1}}$ be the multiplication operators defined by $\Phi_{1,1}$ and $\Phi_{2,1}$ in
\tu{(\ref{REDcoefs})}, respectively, and let $\Ga_{\Phi_{1,2}}$ and $\Ga_{\Phi_{2,2}}$
be the operators associated with $\Phi_{1,2}$ and $\Phi_{2,2}$ in \tu{(\ref{REDcoefs})}
via \tu{(\ref{GaH})}, respectively. Then the operator
\begin{equation}\label{multop}
\mat{cc}{M_{\Phi_{1,1}}&\Ga_{\Phi_{1,2}}\\M_{\Phi_{2,1}}&\Ga_{\Phi_{2,2}}}
:\mat{c}{H^2(\sD_{\om^*})\\\sU}\to\mat{cc}{H^2(\sG)\\ H^2(\sY)}
\end{equation}
is a co-isometry. Moreover, the operator in \tu{(\ref{multop})}
 is unitary whenever $\om$ is an isometry and $\om_2\Pi_\sF$ is
pointwise stable (i.e., $(\om_2\Pi_\sF)^nu\to0$ as $n\to\infty$ for each $u\in\sU$).
\end{theorem}

The addition to the result of Proposition 5.1.6 in \cite{tH07} is that the operator (\ref{multop})
is a co-isometry, rather than just a contraction. To see that this is the case one can follow the
proof of Proposition 5.1.6 in \cite{tH07}, but now also using Theorem \ref{th:ST} below and the
fact that the operator matrix
\begin{equation*}
\mat{cc}{\om\Pi_\sF&D_{\om^*}\\\Pi_\sG&0}:\mat{c}{\sU\\\sD_{\om^*}}\to\mat{c}{\sY\oplus\sU\\\sG}
\end{equation*}
is a co-isometry, which follows from Douglas factorization lemma \cite{D66}.

In order to state Theorem \ref{th:ST} we require some results from system
theory. The terminology corresponds to that in \cite{FFGK98}. A
{\em co-isometric system} is a quadruple $\{A,B,C,D\}$, consisting of
operators $A$ on a Hilbert space $\sX$, $B$ from $\sV$ to $\sX$, $C$
from $\sX$ to $\sW$ and $D$ mapping $\sV$ into $\sW$
such that the operator matrix
\begin{equation}\label{SysMat}
\mat{cc}{A&B\\C&D}:\mat{c}{\sX\\\sV}\to\mat{c}{\sX\\\sW}
\end{equation}
is a co-isometry.
Since $A$ is contractive, we can define functions $F$ and $W$ on $\BD$ by
\begin{equation*}
F(\la)=D+\la C(I-\la A)^{-1}B\ands
W(\la)=C(I-\la A)^{-1}\quad\quad(\la\in\BD).
\end{equation*}
We refer to $F$ and $W$ as the {\em transfer function} and {\em observability function}
associated with $\{A,B,C,D\}$, respectively. {}From the fact that (\ref{SysMat}) is a co-isometry
it follows that $F\in\eS(\sV,\sW)$ and $W\in\eH^2_\tu{ball}(\sX,\sW)$. In particular,
$F$ defines a contractive multiplication operator $M_F$ from $H^2(\sV)$ to $H^2(\sW)$, and
$W$ defines a contraction $\Ga_W$ from $\sX$ into $H^2(\sW)$ via (\ref{GaH}).

\begin{theorem}\label{th:ST}
Let $\{A,B,C,D\}$ be a co-isometric system with transfer function $F\in\eS(\sV,\sW)$
and observability function $W\in\eH^2_\tu{ball}(\sX,\sW)$. Then
\begin{equation}\label{opmat}
\mat{cc}{\!\!M_F&\Ga_W\!\!}:\mat{c}{H^2(\sV)\\\sX}\to H^2(\sW)
\end{equation}
is a co-isometry.
\end{theorem}

\begin{proof}[\bf Proof]
Via the Fourier transform we can identify $H^2(\sV)$ with $\ell^2_+(\sV)$ and $H^2(\sW)$
with $\ell^2_+(\sW)$, where $\ell^2_+(\sV)$ and $\ell^2_+(\sW)$ are the Hilbert spaces
of square summable unilateral sequences of vectors in $\sW$ and $\sV$, respectively.
It then follows that the operator (\ref{opmat}) is a co-isometry if and only if
the operator $\mat{cc}{\!\!T_F&\tilde\Ga_W\!\!}$ from $\ell^2_+(\sV)\oplus\sX$ to $\ell^2_+(\sW)$
is a co-isometry, where $T_F$ and $\tilde\Ga_W$ are given by
\begin{equation*}
T_F=\mat{ccccc}
{D&0&0&0&\cdots\\
CB&D&0&0&\cdots\\
CAB&CB&D&0&\cdots\\
CA^2B&CAB&CB&D&\ddots\\
\vdots&\vdots&\vdots&\ddots&\ddots}
\ands\tilde\Ga_W=\mat{c}{C\\CA\\CA^2\\CA^3\\\vdots}.
\end{equation*}
So we have to show that $X:=T_FT_F^*+\tilde\Ga_W\tilde\Ga_W^*=I$.
It is known that $\mat{cc}{\!\!T_F&\tilde\Ga_W\!\!}$ is contractive (specify Proposition
1.7.2 in \cite{P99} to the time-invariant case), and thus that $X$ is contractive.
Hence it suffices to prove that the $n^\tu{th}$ diagonal entries $X_{n,n}$ of $X$ is equal to $I_\sW$
for each nonnegative integer $n$.
Notice that
\begin{equation*}
X_{n,n}=DD^*+\sum_{i=0}^{n-1}CA^iBB^*A^{*i}C^*+CA^nA^{*n}C^*.
\end{equation*}
Since (\ref{SysMat}) is a co-isometry, we have $BB^*=I-AA^*$ and $DD^*+CC^*=I$. So
\begin{equation*}
\sum_{i=0}^{n-1}CA^iBB^*A^{*i}C^*=CC^*-CA^nA^{*n}C^*.
\end{equation*}
This implies that
\begin{equation*}
X_{n,n}=DD^*+CC^*=I_\sW,
\end{equation*}
which proves our claim.
\end{proof}

{}From the fact that the coefficient matrix (\ref{multop}) is a co-isometry we immediately
obtain the following corollary.

\begin{corollary}\label{cor:nulltocoiso}
Let $\om$ be a contraction of the form \tu{(\ref{givenom})}, and let $\Phi_{2,1}$ and
$\Phi_{2,2}$ be the functions given by \tu{(\ref{REDcoefs})}. Then $\Ga_{\Phi_{2,2}}$ is
a co-isometry if and only if $\Phi_{2,1}(\la)=0$ for each $\la\in\BD$.
\end{corollary}

\section{Existence of a unique solution}\label{sec:UniqueSolution}

In this section we prove Theorem \ref{mainth2}, and give an example that shows that the
earlier obtained sufficient condition ``$\sF=\sU$, $\sY=\{0\}$ or $\om$ is a co-isometry''
is not a necessary condition. We conclude the section with a proposition that shows how
condition (ii) in Theorem \ref{mainth2} simplifies in two special cases.

\begin{proof}[\bf Proof of Theorem \ref{mainth2}]
If there is just one solution, this unique solution must be the central solution
 $H_c$ in (\ref{centralCI}) (or $\Phi_{2.2}$ in the Redheffer representation (\ref{REDsols})).
It thus follows that there is a unique solution if and only if the second summand in (\ref{REDsols})
is the zero function for each $V\in\eS(\sG,\sD_{\om^*})$, that is, if and only if
\begin{equation}\label{secsum}
\Phi_{2,1}(\la)V(\la)(I-\Phi_{1,1}(\la)V(\la))^{-1}\Phi_{1,2}(\la)=0
\quad(\la\in\BD,V\in\eS(\sG,\sD_{\om^*})).
\end{equation}
Hence there is a unique solution whenever
\begin{equation}\label{altcons}
\Phi_{1,2}(\la)=0\quad(\la\in\BD)\qquad\text{or}\qquad\Phi_{2,1}(\la)=0\quad(\la\in\BD).
\end{equation}

The definition of $\Phi_{1,2}$ shows that $\Phi_{1,2}(\la)=0$ for all $\la\in\BD$ if and only if
$\sG=\{0\}$, or equivalently, $\sF=\sU$. Next, from Corollary \ref{cor:nulltocoiso} and
the remark in the paragraph preceding Corollary \ref{cor:coisosolution}
we see that the second sufficient condition in (\ref{altcons})
is equivalent to condition (ii) in Theorem \ref{mainth2}.
This proves that conditions (i) and (ii) in Theorem \ref{mainth2} are sufficient.


It remains to prove the necessity. In fact, from the first part of the proof we see that it
suffices to show that (\ref{secsum}) implies (\ref{altcons}).
Assume that (\ref{secsum}) holds, and assume in addition that $\Phi_{1,2}$ is not the zero function.
Then there exist a $u\in\sU$ and $\la_0\in\BD$ such that $u\not=0$ and $\Phi_{1,2}(\la_0)u\not=0$.
By the continuity of $\Phi_{1,2}$ there exists a $\rho>0$ with $\rho<1-|\la_0|$  such that
$\Phi_{1,2}(\la)u\not=0$ for all $\la$ in the open disc
$\De(\la_0,\rho)=\{\la\in\BD\mid |\la-\la_0|<\rho\}\subset\BD$.

Now fix a $\la\in\De(\la_0,\rho)$ and a $h\in\sD_{\om^*}$. Let $\de>0$ be small enough so that
\begin{equation*}
\|\Phi_{1,2}(\la)u+\de\Phi_{1,1}(\la)h\|\geq\de\|h\|.
\end{equation*}
Then we can define a contraction $U$ from $\sG$ to $\sD_{\om^*}$ by
$U(\Phi_{1,2}(\la)u+\de\Phi_{1,1}(\la)h)=\de h$ and $U k=0$ for each $k\in\sG$ perpendicular to
the vector $\Phi_{1,2}(\la)u+\de\Phi_{1,1}(\la)h$. So
\begin{equation*}
\de(I-U\Phi_{1,1}(\la))h=U\Phi_{1,2}(\la)u.
\end{equation*}
Since $\Phi_{1,1}(\la)$ is a strict contraction (see the second paragraph after Theorem
\ref{th:Redheffer}), we obtain that
\begin{equation*}
U(I-\Phi_{1,1}(\la)U)^{-1}\Phi_{1,2}(\la)u=(I-U\Phi_{1,1}(\la))^{-1}U\Phi_{1,2}(\la)u=\de h.
\end{equation*}
Now let $V\in\eS(\sG,\sD_{\om^*})$ be the constant function with value $U$. Then
\begin{equation*}
\Phi_{2,1}(\la)h=\de^{-1}\Phi_{2,1}(\la)V(\la)(I-\Phi_{1,1}(\la)V(\la))^{-1}\Phi_{1,2}(\la)u=0.
\end{equation*}
Since $h$ is an arbitrary vector in $\sD_{\om^*}$, we have $\Phi_{2,1}(\la)=0$ for each $\la$
in the open disc $\De(\la_0,\rho)$. But $\Phi_{2,1}$ is analytic on $\BD$, thus
$\Phi_{2,1}(\la)=0$ for each $\la\in\BD$.
\end{proof}

\begin{example}\label{ex:counter}
\tu{
Let $\sY=\BC$ and $\sU=\ell^2_+$, that is, $\sU$ is the Hilbert space of square summable unilateral
sequences of complex numbers. Set $\sF=S\ell^2_+$, where $S$ denotes the unilateral forward shift
on $\ell^2_+$. Now define
\begin{equation*}
\om=\mat{cc}{\om_1\\\hline\om_2}
=\mat{cccc}{1&0&0&\cdots\\\hline 0&0&0&\cdots\\0&1&0&\cdots\\0&0&1&\ddots\\\vdots&\vdots&\ddots&\ddots}
:\sF\to\mat{c}{\BC\\\hline\ell^2_+}.
\end{equation*}
Then $\sF\not=\ell^2_+=\sU$, $\om$ is not a co-isometry and $\sY=\BC\not=\{0\}$. However, we do have
\begin{equation*}
\om_1(\Pi_\sF\om_2)^n=\mat{ccccccc}{0&\cdots&0&1&0&0&\cdots}:\sF\to\BC\ \mbox{ for }\ n=0,1,2,\ldots,
\end{equation*}
where the number 1 is situated in the $n^\tu{th}$ position. In particular, $\om_1(\Pi_\sF\om_2)^n$ is
a co-isometry for each nonnegative integer $n$.
So by Theorem \ref{mainth2} the central solution $H_c$ is the only solution to the $H^2$
interpolation problem defined by $\om$.
It is now easy to see that 
$H_c$ is given by
\begin{equation*}
H_c(\la)=\mat{cccccc}{0&1&\la&\la^2&\la^3&\ldots}\quad(\la\in\BD).
\end{equation*}}
\end{example}

\begin{proposition}\label{pr:SpecialCases}
Let $\om$ be a contraction of the form \tu{(\ref{givenom})}. If either $\om_1$ is a strict contraction
or $\dim(\sU)<\infty$, then condition \tu{(ii)} in Theorem \ref{mainth2} is equivalent to $\sY=\{0\}$.
\end{proposition}

\begin{proof}[\bf Proof]
Clearly $\sY=\{0\}$ implies that $\om_1(\Pi_\sF\om_2)^n$ from $\sF$ into $\{0\}$ is a co-isometry.
If $\om_1$ is a strict contraction, then $\om_1$ can only be a co-isometry if $\sY=\{0\}$.
Hence in that case condition (ii) in Theorem \ref{mainth2} reduces to $\sY=\{0\}$.

Next, assume that $\dim(\sU)<\infty$. We already proved that $\sY=\{0\}$ implies condition (ii)
in Theorem \ref{mainth2}. Now assume that condition (ii) in Theorem \ref{mainth2} holds.
Set $\sF_n=\ov{(\om_2^*\Pi_\sF^*)^n\om_1^*\sY}$ for $n=0,1,2,\ldots$. Since
$\om_1(\Pi_\sF\om_2)^n$ is a co-isometry, $\dim(\sF_n)=\dim(\sY)$. Moreover, $\om_1$ maps
$\sF_0$ isometrically onto $\sY$ and $\om_2$ maps $\sF_{n+1}$ isometrically onto $\sF_n$.
Hence the subspaces $\sF_0,\sF_1,\ldots$ of $\sF$ are mutually orthogonal and all of
the same dimension as $\sY$. This is in contradiction with the fact that $\sF\subset\sU$ is
finite dimensional, unless $\sY=\{0\}$.
\end{proof}

\section{Relaxed commutant lifting}\label{sec:RCL}

In this section we return to the relaxed commutant lifting setting. We consider two special cases.
\medskip

\noindent{\bf The sub-optimal case.}
For most of the metric constrained interpolation problems that fit into the commutant
lifting setting there is an interesting special case known as the sub-optimal case;
cf., \cite{F87} for the Nehari problem.
On the level of commutant lifting this corresponds
to the operator $A$ being a strict contraction,
in which case a more explicit description of
all solutions can be obtained (Section XIV.7 in \cite{FF90}, Section VI.6 in \cite{FFGK98}).
To achieve a similar result in the relaxed commutant lifting setting
it is in addition assumed that $R$ has a left inverse;
see \cite{tH}. In that case we have the following result.

\begin{corollary}\label{cor:suboptimal}
Let $\LDS$ be a data set with $A$ a strict contraction and $R$ left invertible.
Then there exists a unique contractive interpolant for $\LDS$
if and only if $\ov{Q\sH_0}=\sH$  or $T'$ is an isometry.
\end{corollary}

\begin{proof}[\bf Proof]
Note that $A$ being a strict contraction corresponds to $D_A$ being invertible on $\sH$. The second
condition in (\ref{intertw}) implies that $Q$ has a left inverse whenever $R$ has a left inverse.
So the conditions on the data set imply that both $D_AQ$ and $D_AR$ have left inverses, and thus
\begin{equation*}
\om_2=D_AR(Q^*D_A^2Q)^{-1}Q^*D_A|_\sF
\end{equation*}
has a left inverse as well. The latter implies that $\om_1$ is a strict
contraction, and thus, by Proposition \ref{pr:SpecialCases} condition (ii) in
Theorem \ref{mainth1} reduces to $\sD_{T'}=\{0\}$, that is, $T'$ is an isometry.
Since $D_A$ is invertible, it follows that condition (i) in Theorem \ref{mainth1}
is equivalent to $\ov{Q\sH_0}=\sH$.
\end{proof}

\noindent{\bf Relaxations of metric constrained interpolation problems.}
As a motivation for the relaxed commutant lifting problem, in \cite{FFK02} a number of classical
metric constrained interpolation problems, including Nevanlinna-Pick and Sarason interpolation,
were provided with a relaxed version.
The common ingredient in all these problems is that, in the relaxed commutant lifting setting, $R$ and $Q$
are operators from $\sV^{n-1}$ into $\sV^n$ (where $\sV$ is some Hilbert space and $n$ a positive
integer) of the form
\begin{equation}\label{RandQ}
R=\mat{cccc}
{I_\sV&0&\cdots&0\\
0&I_\sV&\ddots&\vdots\\
\vdots&\ddots&\ddots&0\\
0&\cdots&0&I_\sV\\
0&\cdots&0&0}\ands
Q=\mat{cccc}
{0&0&\cdots&0\\
I_\sV&0&\cdots&0\\
0&\ddots&\ddots&\vdots\\
\vdots&\ddots&I_\sV&0\\
0&\cdots&0&I_\sV}.
\end{equation}
In particular, if $\dim(\sV)<\infty$, then $\dim(\sH)=\dim(\sV^n)=\dim(\sV)^n<\infty$, and thus
$\dim(\sD_A)<\infty$. So it follows from Proposition \ref{pr:SpecialCases} that
in this case condition (ii) of Theorem \ref{mainth1} reduces to $\sD_{T'}=\{0\}$, i.e., $T'$
is an isometry.

A typical result for the relaxed interpolation problems of \cite{FFK02} is that in the special case
$\sV=\BC$ there exists a unique contractive interpolant whenever the operator $A$  has norm one.
We show here that this is the case for each relaxed commutant lifting problem with $R$ and $Q$ are
of the form (\ref{RandQ}) and $\sV=\BC$. First we prove the following lemma.

\begin{lemma}\label{lem:refcon1}
Let $\LDS$ be a data set for the relaxed commutant lifting problem, and set $\sG=\sD_A\ominus\sF$.
Then $\ov{D_A\sG}$ is perpendicular to both $Q\sH_0$ and $\kr D_A$. In particular, $\sF=\sD_A$ whenever
 $\ov{Q\sH_0}\vee\kr D_A=\sH$.
\end{lemma}

\begin{proof} Let $g\in\sG=\kr Q^*D_A\subset\sD_A$. Then $D_A g\in\kr Q^*$, or equivalently,
$D_Ag$ is perpendicular to $Q\sH_0$. Since $D_A$ is a self adjoint operator, it follows
directly that $D_Ag$ is perpendicular to $\kr D_A$. The last statement of Lemma \ref{lem:refcon1}
holds because $D_A$ restricted to $\sG$ is one-to-one.
\end{proof}

\begin{corollary}\label{cor:specialRQ}
Let $\LDS$ be a data set for the relaxed commutant lifting problem with $R$ and $Q$ given by
\tu{(\ref{RandQ})}, $\sV=\BC$ and $\sD_{T'}\not=\{0\}$. Then there exists a unique contractive
interpolant if and only if $A$ has norm one.
\end{corollary}

\begin{proof} In case the norm of $A$ is less than one, we are in the sub-optimal case described
above, while $\ker Q^*\not=\{0\}$ and $\sD_{T'}\not=\{0\}$. Thus there is more than one contractive
interpolant.

Now assume that $\|A\|=1$. Since $\sH=\BC^n$ is finite dimensional, this implies that $A$ has a
norm attaining vector, i.e., $\kr D_A\not=\{0\}$. Note that
$\kr Q^*=\BC\oplus\{0\}^{n-1}$.
So according to Lemma \ref{lem:refcon1} it suffices to show that there exists a norm attaining
vector $h=(h_1,\ldots,h_n)\in\BC^n$ for $A$ with $h_1\not=0$.

Let $h=(h_1,\ldots,h_n)\not=0$ be a norm attaining vector, but assume that $h_1=0$. Then
$h=Qk$, with $k=(h_2,\ldots,h_n)\in\BC^{n-1}$. We have
\begin{equation*}
\|Rk\|=\|h\|=\|Ah\|=\|AQk\|=\|T'ARk\|\leq\|ARk\|\leq\|Rk\|.
\end{equation*}
Thus $Rk=(h_2,\ldots,h_n,0)$ is also a norm attaining vector. It may happen
that $h_2=0$, in which case we just repeat the above procedure. After at most $n-1$ times we arrive
at a norm attaining vector $h=(h_1,\ldots,h_n)$ with $h_1\not=0$.
\end{proof}

\noindent\textbf{Acknowledgement.}
The author thanks Prof.~M.A.
Kaashoek and Prof.~A.E. Frazho for their useful comments and
suggestions.

\end{document}